**Predicting some physicochemical properties of octane isomers: A topological approach using *ev*-degree and *ve*-degree Zagreb indices**

Süleyman Ediz[1]

**Abstract** Topological indices have important role in theoretical chemistry for QSPR researches. Among the all topological indices the Randić and the Zagreb indices have been used more considerably than any other topological indices in chemical and mathematical literature. Most of the topological indices as in the Randić and the Zagreb indices are based on the degrees of the vertices of a connected graph. Recently novel two degree concepts have been defined in graph theory; *ev*-degrees and *ve*-degrees. In this study we define *ev*-degree Zagreb index, *ve*-degree Zagreb indices and *ve*-degree Randić index by using these new graph invariants as parallel to their corresponding classical degree versions. We compare these new group *ev*-degree and *ve*-degree indices with the other well-known and most used topological indices in literature such as; Wiener, Zagreb and Randić indices by modelling some physicochemical properties of octane isomers. We show that the *ev*-degree Zagreb index, the *ve*-degree Zagreb and the *ve*-degree Randić indices give better correlation than Wiener, Zagreb and Randić indices to predict the some specific physicochemical properties of octanes. We investigate the relations between the second Zagreb index and *ev*-degree and *ve*-degree Zagreb indices and some mathematical properties of *ev*-degree and *ve*-degree Zagreb indices.

# 1 Introduction

Graph theory which is an important branch of applied mathematics has many applications to modelling real world problems from science to technology. Chemical graph theory which is a fascinating branch of graph theory has many applications related to chemistry. Chemical graph theory provides many information about molecules and atoms by using pictorial representation (chemical graph) of these chemical compounds. A topological index which is a numerical quantity derived from the chemical graph of a molecule is used to modelling chemical and physical properties of molecules in QSPR/QSAR researches. Quantitative structure-property/activity relationships (QSPR/QSAR) studies have very important role in theoretical chemistry. Octane isomers have been used widely in QSPR studies. The role of octane isomers in QSPR studies, we refer the interested reader [1-4] and references therein. Among the all topological indices, Wiener, Randić and Zagreb indices are the most used topological indices in the chemical and mathematical literature so far.

[1] Faculty of Education, Yuzuncu Yil University, Van, Turkey
e-mail: suleymanediz@yyu.edu.tr

Very recently, Chellali, Haynes, Hedetniemi and Lewis have published a seminal study: On *ve*-degrees and *ev*-degrees in graphs [5]. The authors defined two novel degree concepts in graph theory; *ev*-degrees and *ve*-degrees and investigate some basic mathematical properties of both novel graph invariants with regard to graph regularity and irregularity [5]. After given the equality of the total *ev*-degree and total *ve*-degree for any graph, also the total *ev*-degree and the total *ve*-degree were stated as in relation to the first Zagreb index. It was proposed in the article that the chemical applicability of the total *ev*-degree (and the total *ve*-degree) could be an interesting problem in view of chemistry and chemical graph theory.

In this study we define *ev*-degree Zagreb index, *ve*-degree Zagreb indices and *ve*-degree Randić index by using these new graph invariants. We define these novel topological indices invariants as parallel to corresponding original definitions of based on classical degree concept. We compare these new group *ev*-degree, *ve*-degree Zagreb and *ve*-degree Randić indices with the other well-known and most used topological indices such as Wiener, Zagreb and Randić indices by modelling some physicochemical properties of octane isomers.

## 2 Preliminaries

In this section we give some basic and preliminary concepts which we shall use later. A graph $G = (V, E)$ consists of two nonempty sets $V$ and 2-element subsets of $V$ namely $E$. The elements of $V$ are called vertices and the elements of $E$ are called edges. For a vertex $v$, $\deg(v)$ show the number of edges that incident to $v$. The set of all vertices which adjacent to $v$ is called the open neighborhood of $v$ and denoted by $N(v)$. If we add the vertex $v$ to $N(v)$, then we get the closed neighborhood of $v$, $N[v]$. For the vertices $u$ and $v$, $d(u, v)$ denotes the distance between $u$ and $v$ which means that minimum number of edges between $u$ and $v$. In [6], the Wiener index of a connected graph $G$, the first topological index, was defined as;

$$W = W(G) = \frac{1}{2}\sum_{u,v \in V(G)} d(u,v).$$

In his study, Wiener used the total distance between all different atoms (vertices) of paraffin to predict boiling point. We refer the interested reader to [7-9] and the references therein for the detailed discussion of Wiener index.

The first and second Zagreb indices [10] defined as follows: The first Zagreb index of a connected graph $G$, defined as;

$$M_1 = M_1(G) = \sum_{u \in V(G)} \deg(u)^2 = \sum_{uv \in E(G)} (\deg(u) + \deg(v)).$$

And the second Zagreb index of a connected graph $G$, defined as;

$$M_2 = M_2(G) = \sum_{uv \in E(G)} \deg(u) \cdot \deg(v).$$

The authors investigated the relationship between the total $\pi$-electron energy on molecules and Zagreb indices [10]. For the details see the references [11-13]. Randić investigated the measuring the extent of branching of the carbon-atom skeleton of saturated hydrocarbons via Randić index [14]. The Randić index of a connected graph $G$ defined as;

$$R = R(G) = \sum_{uv \in E(G)} (\deg(u) \cdot \deg(v))^{-1/2}.$$

We refer the interested reader to [15-17] and the references therein for the up to date arguments about the Randić index. And now we give the definitions of *ev*-degree and *ve*-degree concepts which were given by Chellali *et al.* in [5].

**Definition 2.1** [5] Let $G$ be a connected graph and $v \in V(G)$. The *ve*-degree of the vertex $v$, $deg_{ve}(v)$, equals the number of different edges that incident to any vertex from the closed neighborhood of $v$. For convenience we prefer to show the *ve*-degree of the vertex $v$, $c_v$.

**Definition 2.2** [5] Let $G$ be a connected graph and $e = uv \in E(G)$. The *ev*-degree of the edge $e$, $deg_{ev}(e)$, equals the number of vertices of the union of the closed neighborhoods of $u$ and $v$. For convenience we prefer to show the *ev*-degree of the edge $e = uv$, $c_e$ or $c_{uv}$.

We illustrate these new degree definitions for the vertices and edges of the graph $G$ which are shown in Figure 1.

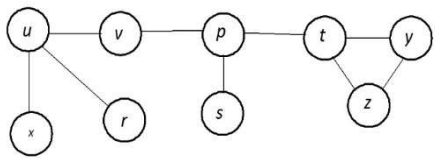

**Figure 1** The graph *G* for the Example 2.3 and Example 2.13

**Example 2.3** Notice that for the vertices of $G$, we get $c_x = 3$, $c_u = 4$, $c_v = 6$, $c_r = 3$, $c_p = 6$, $c_s = 3$, $c_t = 6$, $c_z = 4$ and $c_y = 4$. And for the edges of $G$, we get $c_{xu} = 4$, $c_{uv} = 5$, $c_{ur} = 4$, $c_{vp} = 5$, $c_{pt} = 6$, $c_{ps} = 4$, $c_{tz} = 4$, $c_{ty} = 4$ and $c_{yz} = 3$.

**Definition 2.4** [5] Let $G$ be a connected graph and $v \in V(G)$. The total *ev*-degree of the graph $G$ is defined as;

$$T_e = T_e(G) = \sum_{e \in E(G)} c_e.$$

And the total *ve*-degree of the graph $G$ is defined as;

$$T_v = T_v(G) = \sum_{v \in V(G)} c_v.$$

**Observation 2.5** [5] For any connected graph G,

$$T_e(G) = T_v(G).$$

The following theorem states the relationship between the first Zagreb index and the total *ve*-degree of a connected graph $G$.

**Theorem 2.6** [5] For any connected graph G,

$$T_e(G) = T_v(G) = M_1(G) - 3n(G).$$

where $n(G)$ denotes the total number of triangles in G.

We can restate the Theorem 2.1 for the trees which are acyclic and are not contain any triangles.

**Corollary 2.7** For any tree $T$,

$$T_e(T) = T_v(T) = M_1(T).$$

And from this last equality we naturally consider to apply these two novel degree concepts to chemical graph theory by introducing *ev*-degree and *ve*-degree Zagreb indices as well as *ve*-degree Randić index.

**Definition 2.8** Let $G$ be a connected graph and $e \in E(G)$. The *ev*-degree Zagreb index of the graph $G$ is defined as;

$$S = S(G) = \sum_{e \in E(G)} c_e^2.$$

**Definition 2.9** Let $G$ be a connected graph and $v \in V(G)$. The first *ve*-degree Zagreb alpha index of the graph $G$ is defined as;

$$S^\alpha = S^\alpha(G) = \sum_{v \in V(G)} c_v^2.$$

**Definition 2.10** Let $G$ be a connected graph and $uv \in E(G)$. The first *ve*-degree Zagreb beta index of the graph $G$ is defined as;

$$S^\beta = S^\beta(G) = \sum_{uv \in E(G)} (c_u + c_v).$$

**Definition 2.11** Let $G$ be a connected graph and $uv \in E(G)$. The second *ve*-degree Zagreb index of the graph $G$ is defined as;

$$S^\mu = S^\mu(G) = \sum_{uv \in E(G)} c_u c_v.$$

**Definition 2.12** Let $G$ be a connected graph and $uv \in E(G)$. The *ve*-degree Randić index of the graph $G$ is defined as;

$$R^\alpha(G) = \sum_{uv \in E(G)} (c_u c_v)^{-1/2}.$$

**Example 2.13** We compute these novel topological indices for the graph $G$ in the Example 2.3 (see Figure 2.1).

$S = S(G) = \sum_{e \in E(G)} c_e^2 = 175$, $S^\alpha = S^\alpha(G) = \sum_{v \in V(G)} c_v^2 = 183$, $S^\beta = S^\beta(G) = \sum_{uv \in E(G)} (c_u + c_v) = 84$,

$S^\mu = S^\mu(G) = \sum_{uv \in E(G)} c_u c_v = 202$, $R^\alpha(G) = \sum_{uv \in E(G)} 1/(c_u c_v)^{1/2} = 13.425$ and $M_2(G) = 46$.

## 3 Results and Discussions

In this section we compare all above mentioned old and new topological indices with each other by using strong correlation coefficients acquired from the chemical graphs of octane isomers. We get the experimental results at the www.moleculardescriptors.eu (see Table 1). The following physicochemical features have been modeled:

• Entropy,

• Acentric factor (AcenFac),

• Enthalpy of vaporization (HVAP),

• Standard enthalpy of vaporization (DHVAP).

We select those physicochemical properties of octane isomers for which give reasonably good correlations, i.e. the absolute value of correlation coefficients are larger than 0.8 except from the property HVAP (see Table 2). Also we find the Wiener index, the first Zagreb index, the second Zagreb index and the Randić indices of octane isomers values at the www.moleculardescriptors.eu (see Table 3). We also calculate and show the *ev*-degree Zagreb index, the *ve*-degree Zagreb indices and the *ve*-degree Randić index of octane isomers values in Table 3.

Table 1. Some physicochemical properties of octane isomers

| Molecule | Entropy | AcenFac | HVAP | DHVAP |
|---|---|---|---|---|
| n-octane | 111.70 | 0.39790 | 73.19 | 9.915 |
| 2-methyl-heptane | 109.80 | 0.37792 | 70.30 | 9.484 |
| 3-methyl-heptane | 111.30 | 0.37100 | 71.30 | 9.521 |
| 4-methyl-heptane | 109.30 | 0.37150 | 70.91 | 9.483 |
| 3-ethyl-hexane | 109.40 | 0.36247 | 71.70 | 9.476 |
| 2,2-dimethyl-hexane | 103.40 | 0.33943 | 67.70 | 8.915 |
| 2,3-dimethyl-hexane | 108.00 | 0.34825 | 70.20 | 9.272 |
| 2,4-dimethyl-hexane | 107.00 | 0.34422 | 68.50 | 9.029 |
| 2,5-dimethyl-hexane | 105.70 | 0.35683 | 68.60 | 9.051 |
| 3,3-dimethyl-hexane | 104.70 | 0.32260 | 68.50 | 8.973 |
| 3,4-dimethyl-hexane | 106.60 | 0.34035 | 70.20 | 9.316 |
| 2-methyl-3-ethyl-pentane | 106.10 | 0.33243 | 69.70 | 9.209 |
| 3-methyl-3-ethyl-pentane | 101.50 | 0.30690 | 69.30 | 9.081 |
| 2,2,3-trimethyl-pentane | 101.30 | 0.30082 | 67.30 | 8.826 |
| 2,2,4-trimethyl-pentane | 104.10 | 0.30537 | 64.87 | 8.402 |
| 2,3,3-trimethyl-pentane | 102.10 | 0.29318 | 68.10 | 8.897 |
| 2,3,4-trimethyl-pentane | 102.40 | 0.31742 | 68.37 | 9.014 |
| 2,2,3,3-tetramethylbutane | 93.06 | 0.25529 | 66.20 | 8.410 |

Table 2. The correlation coefficients between new and old topological indices and some physicochemical properties of octane isomers

| Index | Entropy | AcenFac | HVAP | DHVAP |
|---|---|---|---|---|
| S | -0.9614 | -0.9829 | -0.8425 | -0.9043 |
| $S^\alpha$ | -0.9565 | -0.9906 | -0.8279 | -0.8931 |
| $S^\beta$ | -0.9410 | -0.9864 | -0.7281 | -0.8118 |
| $S^\mu$ | -0.9481 | -0.9863 | -0.7552 | -0.8118 |
| $R^\alpha$ | 0.9486 | 0.9829 | 0.8351 | 0.8924 |
| W | 0.8772 | 0.9656 | 0.7381 | 0.8202 |
| $M_1$ | -0.9543 | -0.9731 | -0,8860 | -0.9361 |
| $M_2$ | -0.9410 | -0.9864 | -0.7281 | -0.8118 |
| R | 0.9063 | 0.9043 | 0.9359 | 0.9580 |

**Table 3.** Topological indices of octane isomers

| Molecule | $M_1$ | $M_2$ | W | R | S | $S^\alpha$ | $S^\beta$ | $S^\mu$ | $R^\alpha$ |
|---|---|---|---|---|---|---|---|---|---|
| n-octane | 26 | 24 | 84 | 3.914 | 98 | 90 | 48 | 84 | 2.144 |
| 2-methyl-heptane | 28 | 26 | 79 | 3.770 | 114 | 104 | 52 | 98 | 1.971 |
| 3-methyl-heptane | 28 | 27 | 76 | 3.808 | 116 | 98 | 54 | 106 | 1.956 |
| 4-methyl-heptane | 28 | 27 | 75 | 3.808 | 116 | 110 | 54 | 107 | 1.991 |
| 3-ethyl-hexane | 28 | 28 | 72 | 3.846 | 118 | 114 | 56 | 115 | 1.964 |
| 2,2-dimethyl-hexane | 32 | 30 | 71 | 3.561 | 152 | 138 | 60 | 132 | 1.754 |
| 2,3-dimethyl-hexane | 30 | 30 | 70 | 3.681 | 134 | 126 | 60 | 129 | 1.784 |
| 2,4-dimethyl-hexane | 30 | 29 | 71 | 3.664 | 132 | 124 | 58 | 121 | 1.799 |
| 2,5-dimethyl-hexane | 30 | 28 | 74 | 3.626 | 130 | 118 | 56 | 113 | 1.801 |
| 3,3-dimethyl-hexane | 32 | 32 | 67 | 3.621 | 156 | 146 | 64 | 148 | 1.718 |
| 3,4-dimethyl-hexane | 30 | 31 | 68 | 3.719 | 136 | 130 | 62 | 136 | 1.753 |
| 2-methyl-3-ethyl-pentane | 30 | 31 | 67 | 3.719 | 136 | 132 | 62 | 137 | 1.770 |
| 3-methyl-3-ethyl-pentane | 32 | 34 | 64 | 3.682 | 160 | 152 | 68 | 163 | 1.645 |
| 2,2,3-trimethyl-pentane | 34 | 35 | 63 | 3.481 | 174 | 162 | 70 | 171 | 1.527 |
| 2,2,4-trimethyl-pentane | 34 | 32 | 66 | 3.417 | 168 | 156 | 64 | 147 | 1.606 |
| 2,3,3-trimethyl-pentane | 34 | 36 | 62 | 3.504 | 176 | 164 | 72 | 179 | 1.489 |
| 2,3,4-trimethyl-pentane | 32 | 33 | 65 | 3.553 | 152 | 144 | 66 | 151 | 1.589 |
| 2,2,3,3-tetramethylbutane | 38 | 40 | 58 | 3.250 | 214 | 194 | 80 | 217 | 1.277 |

**Table 4.** The squares of correlation coefficients between topological indices and some physicochemical properties of octane isomers

| Index | Entropy | AcenFac | HVAP | DHVAP |
|---|---|---|---|---|
| S | 0.9242 | 0.9660 | 0.7098 | 0.8177 |
| $S^\alpha$ | 0.9148 | 0.9812 | 0.6854 | 0.7976 |
| $S^\beta$ | 0.8854 | 0.9729 | 0.5301 | 0.6590 |
| $S^\mu$ | 0.8988 | 0.9727 | 0.5703 | 0.6590 |
| $R^\alpha$ | 0.8998 | 0.9660 | 0.6973 | 0.7963 |
| W | 0.7694 | 0.9323 | 0.5447 | 0.6727 |
| $M_1$ | 0.9106 | 0.9469 | 0.7849 | 0.8762 |
| $M_2$ | 0.8854 | 0.9729 | 0.5301 | 0.6590 |
| R | 0.8213 | 0.8177 | 0.8759 | 0.9177 |

It can be seen from the Table 2 that the most convenient indices which are modelling the Entropy, Enthalpy of vaporization (HVAP), Standard enthalpy of vaporization (DHVAP) and Acentric factor (AcenFac) are *ve*-degree Zagreb index ($S$) for entropy, the first *ve*-degree Zagreb alpha index ($S^\alpha$) for Acentric Factor and the Randić index ($R$) for the Enthalpy of vaporization (HVAP) and Standard enthalpy of vaporization (DHVAP), respectively. But notice that the first two indices show the negative strong correlation and the third index show the positive strong correlation. Because of this fact we compare these graph invariants with each other by using the squares of correlation coefficients for ensure the compliance between the positive and negative correlation coefficients (see Table 4).

**Entropy:** We can see from the Table 4 that the *ve*-degree Zagreb index ($S$) gives the highest square of correlation coefficient for entropy. After that the first *ve*-degree Zagreb alpha index ($S^\alpha$), the first Zagreb index ($M_1$), the *ve*-degree Randić index ($R^\alpha$) and the second *ve*-degree Zagreb index ($S^\mu$) give the highest square of correlation coefficients, respectively.

**Acentric factor (AcenFac):** We can see from the Table 4 that the first *ve*-degree Zagreb alpha index ($S^\alpha$) gives the highest square of correlation coefficient for Acentric factor. After that the first *ve*-degree Zagreb beta index ($S^\beta$) and the second Zagreb index ($M_2$) give the same value. And the the second *ve*-degree Zagreb index ($S^\mu$), the *ev*-degree Zagreb index ($S$) and *ve*-degree Randić index ($R^\alpha$) give the highest square of correlation coefficients, respectively.

**Enthalpy of vaporization (HVAP):** It can be seen from the Table 4 that the Randić index ($R$) gives the the highest square of correlation coefficient for Enthalpy of vaporization. After that the first Zagreb index ($M_1$), the *ev*-degree Zagreb index ($S$), the *ve*-degree Randić index ($R^\alpha$) and the the first *ve*-degree Zagreb alpha index ($S^\alpha$) give the highest square of correlation coefficients, respectively.

**Standard enthalpy of vaporization (DHVAP):** We can observe from the Table 4 that the Randić index ($R$) gives the the highest square of correlation coefficient for Enthalpy of vaporization. After that the first Zagreb index ($M_1$), the *ev*-degree Zagreb index ($S$), the first *ve*-degree Zagreb alpha index ($S^\alpha$) and the *ve*-degree Randić index ($R^\alpha$) give the highest square of correlation coefficients, respectively.

And now we investigate the relations between the old topological indices and the novel topological indices. The correlation coefficients between the Wiener, Zagreb, Randić indices and the ev-degree and ve-degree indices are shown in Table 5. It can be shown from the Table 5 that the first *ve*-degree Zagreb beta index ($S^\beta$) gives the highest absolute value of correlation coefficient with the Wiener index. The *ev*-degree Zagreb index ($S$) gives the highest correlation coefficient with the first Zagreb index ($M_1$). The first *ve*-degree Zagreb beta index ($S^\beta$) gives the highest absolute value of correlation coefficient with the Randić index. And it is very surprisingly see that the correlation coefficient between the second Zagreb index ($M_2$) and the first *ve*-degree Zagreb beta index ($S^\beta$) is one. We can see from the Table 3 that $S^\beta(G) = 2M_2(G)$ for the molecular graphs of octane isomers. But we know that $S^\beta(G) \neq 2M_2(G)$ from the Example 2.3. The following section we investigate the relation between the second Zagreb index and the first *ve*-degree Zagreb beta index.

Table 5. The correlation coefficients between old and corresponding novel topological indices

| Index | W | $M_1$ | $M_2$ | R |
|---|---|---|---|---|
| S | -0.9177 | 0.9951 | 0.9676 | -0.9441 |
| $S^\alpha$ | 0.9483 | 0.9818 | 0.9774 | -0.9182 |
| $S^\beta$ | -0.9683 | 0.9495 | 1.000 | -0.8609 |
| $S^\mu$ | -0.9567 | 0.9523 | 0.9982 | -0.8645 |
| $R^\alpha$ | 0.9478 | -0.9764 | -0.9758 | 0.9365 |

The cross correlation matrix of *ev*-degree and *ve*-degree indices are given in Table 6.

Table 6. The cross correlation matrix of the *ev*-degree and *ve*-degree topological indices

| Index | S | $S^\alpha$ | $S^\beta$ | $S^\mu$ | $R^\alpha$ |
|---|---|---|---|---|---|
| S | 1.0000 | | | | |
| $S^\alpha$ | 0.9901 | 1.0000 | | | |
| $S^\beta$ | 0.9676 | 0.9774 | 1.0000 | | |
| $S^\mu$ | 0.9738 | 0.9797 | 0.9982 | 1.0000 | |
| $R^\alpha$ | -0.9758 | -0.9752 | -0.9758 | -0.9701 | 1.0000 |

It can be shown from the Table 6 that the minimum correlation efficient among the all *ve*-degree and *ev*-degree indices is 0.9676 which is indicate strong correlation among all these novel indices. From the above arguments, we can say that the *ve*-degree and *ev*-degree indices are possible tools for QSPR researches.

**4 Lower and upper bounds of *ev*-degree and *ve*-degree Zagreb indices for general graphs**

In this section are given the relations between second Zagreb index and *ve*-degree and *ev*-degree Zagreb indices. And also fundamental mathematical properties of *ev*-degree and *ve*-degree Zagreb indices are given.

**Lemma 4.1** *Let T be a tree and $v \in V(T)$ then,*

$$c_v = \sum_{u \in N(v)} \deg(u).$$

*Proof* From the Definition 2.1 we know that $c_v$ equals the number of different edges incident to any vertex from $N(v)$. Clearly for any tree, this definition corresponds the sum of all degrees of the vertices lie in $N(v)$. Hence

$$c_v = \sum_{u \in N(v)} \deg(u). \qquad \square$$

**Theorem 4.2** *Let T be a tree with the vertex set $V(T) = \{v_1, v_2, \ldots, v_n\}$ then*

$$S^\beta(T) = 2M_2(T).$$

*Proof* From the Definition 2.10 and Lemma 4.1 we can directly write

$$S^\beta(T) = \sum_{v_i v_j \in E(T)} (c_{v_i} + c_{v_j}) = \sum_{v_i v_j \in E(T)} \left( \sum_{w \in N(v_i)} \deg(w) + \sum_{w \in N(v_j)} \deg(w) \right)$$

$$= \deg(v_1) \sum_{w \in N(v_1)} \deg(w) + \deg(v_2) \sum_{w \in N(v_2)} \deg(w) + \cdots + \deg(v_n) \sum_{w \in N(v_n)} \deg(w)$$

Notice that the above sum contains the multiplication of the degree of end vertices of each edge exactly two times. Hence,

$$= 2 \sum_{uv \in E(T)} \deg(u) \deg(v) = 2M_2(T). \qquad \square$$

Before we give the following interesting theorem, we mention the forgotten topological index [10]. The forgotten topological index for a connected graph $G$ defined as;

$$F = F(G) = \sum_{v \in V(G)} \deg(v)^3 = \sum_{uv \in E(G)} (\deg(u)^2 + \deg(v)^2).$$

It was showed in [18] that the predictive power of the forgotten topological index is very close to the first Zagreb index for the entropy and acentric factor. For further studies about the forgotten topological index we refer to the interested reader [18-20] and references therein.

**Theorem 4.3** *Let G be a triangle free connected graph, then;*

$$S(G) = F(G) + 2M_2(G).$$

*Proof.* It was showed in [5] that $c_e = c_{uv} = \deg(u) + \deg(v)$ for any triangle free graph. By using this equality, we get that;

$$S = S(G) = \sum_{e=uv \in E(G)} c_e^2 = \sum_{e=uv \in E(G)} (\deg(u) + \deg(v))^2$$

$$= \sum_{e=uv \in E(G)} (\deg(u)^2 + \deg(v)^2) + 2 \sum_{e=uv \in E(G)} \deg(u)\deg(v)$$

$$= F(G) + 2M_2(G). \qquad \square$$

We can state the following corollary which describe the relation between the *ev*-degree Zagreb index and the first *ve*-degree Zagreb alpha index for trees by using the Theorem 4.3.

**Corollary 4.4** *Let T be a tree then;*

$$S(T) = F(T) + S^\beta(T).$$

And now we give the maximum and minimum graph classes with respect to *ev*-degree and *ve*-degree Zagreb indices.

**Theorem 4.5** *Let G be a simple connected graph of order $n \geq 3$ vertices then;*

$$16n - 30 \leq S(G) \leq \frac{1}{2}n^3(n-1).$$

*Lower bound is achieved if and only if G is a path and upper bound is achieved if and only if G is a complete graph.*

*Proof* We get that $c_e = c_{uv} = |N(u) \cup N(V)|$ from the definition of *ev*-degree of any edge of $G$. $|N(u) \cup N(V)|$ reaches its maximum value for the complete graphs and its minimum value for the path for an edge of $G$. There are $n - 3$ edges with their *ev*-degrees equals 4 and 2 edges with their *ev*-degrees equals 3 for the *n*-vertex path.

And the *ev*-degrees of all edges of the complete graph are $n$. From this, the desired result is acquired.
□

**Theorem 4.6** *Let T be a tree of order $n \geq 3$ vertices then;*

$$16n - 30 \leq S(T) \leq n^2(n-1).$$

*Lower bound is achieved if and only if T is a path and upper bound is achieved if and only if T is a star.*

*Proof* The lower bound comes from Theorem 4.5. From the same arguments of the Theorem 4.5, the maximum tree of the *ev*-degree Zagreb index is star graph. The *ev*-degrees of all edges of the star graph are $n$. From this, the proof is completed. □

**Theorem 4.7** *Let G be a simple connected graph of order $n \geq 5$ vertices then;*

$$16n - 6 \leq S^\alpha(G) \leq \frac{1}{4}n^3(n-1)^2.$$

*Lower bound is achieved if and only if G is a path and upper bound is achieved if and only if G is a complete graphs.*

*Proof* We know that $c_u$ equals the number of different edges that incident to any vertex from the closed neighborhood of $v$. Clearly $c_u$ reaches its maximum value for the complete graphs and its minimum value for the path for a vertex of $G$. There are $n - 2$ vertices with their *ve*-degrees equals 4, 2 vertices with their *ve*-degrees equals 3 and 2 vertices with their *ve*-degrees equals 2. And the *ve*-degrees of all vertices of the complete graph are $n(n-1)/2$. From this, the desired result is acquired. □

**Theorem 4.8** *Let T be a tree of order $n \geq 5$ vertices then;*

$$16n - 6 \leq S^\alpha(T) \leq n(n-1)^2.$$

*Lower bound is achieved if and only if T is a path and upper bound is achieved if and only if T is a star.*

*Proof* The lower bound comes from Theorem 4.6. From the same arguments of the Theorem 4.7, the maximum tree of the *ve*-degree first Zagreb alpha index is star graph. The *ve*-degrees of all vertices of the star graph are $n - 1$. From this, the proof is completed. □

**Theorem 4.9** *Let G be a simple connected graph of order $n \geq 5$ vertices then;*

$$8n - 16 \leq S^\beta(G) \leq \frac{1}{2}n^2(n-1)^2.$$

*Lower bound is achieved if and only if G is a path and upper bound is achieved if and only if G is a complete graph.*

*Proof* The proof is similar the proof of Theorem 4.7. □

**Theorem 4.10** *Let T be a tree of order $n \geq 5$ vertices then;*

$$16n - 6 \leq S^\beta(G) \leq 2n(n-1).$$

*Lower bound is achieved if and only if T is a path and upper bound is achieved if and only if T is a star.*

*Proof* The proof is similar the proof of Theorem 4.8. □

**Theorem 4.11** *Let G be a simple connected graph of order $n \geq 5$ vertices then;*

$$16n - 44 \leq S^\mu(G) \leq \frac{1}{8}n^3(n-1)^3.$$

*Lower bound is achieved if and only if G is a path and upper bound is achieved if and only if G is a complete graph.*

*Proof* The proof is similar the proof of Theorem 4.7. □

**Theorem 4.12** *Let T be a tree of order $n \geq 5$ vertices then;*

$$16n - 6 \leq S^\mu(T) \leq (n-1)^3.$$

*Lower bound is achieved if and only if T is a path and upper bound is achieved if and only if T is a star.*

*Proof* The proof is similar the proof of Theorem 4.8. □

**5 Conclusion**

We proposed novel topological indices based on *ev*-degree and *ve*-degree concept which have been defined very recently in graph theory. It has been shown that these indices can be used as predictive means in QSAR researches. Predictive power of these indices have been tested on by using some physicochemical properties of octanes. Acquired results show that the new *ev*-degree and *ve*-degree indices give somewhat better results by analogy well-known Wiener, Zagreb and Randić indices. In addition, we investigated basic mathematical properties of these novel topological indices. We have found a lower and upper bounds for the simple connected graphs. It can be interesting to find the exact value of the *ev*-degree and *ve*-degree Zagreb indices of some graph operations such as; direct, Cartesian, corona, tensor, hierarchical and generalized hierarchical product of graphs for further studies. It can also be interesting to investigate the relations between the *ev*-degree and *ve*-degree Zagreb indices and the other well-known topological indices.

**References**

1. M. Kuanar, S.K. Kuanar, B.K. Mishra, I. Gutman, Indian Journal of Chemistry-Section A 38A, 525 (1999).
2. M. Randić, New Journal of Chemistry 20, 1001 (1996).
3. M. Randić, M. Pompe, SAR and QSAR in Environmental Research 10, 451 (1999).
4. M.H. Sunilkumar, Applied Mathematics and Nonlinear Sciences 1,345 (2016).
5. M. Chellali, T.W. Haynes, S.T. Hedetniemi, T.M. Lewis, Discrete Mathematics 340, 31 (2017).
6. H. Wiener, J. Am. Chem. Soc. 69, 17 (1947).
7. R. Liu, X. Du, H. Jia, Bull. Aust. Math. Soc. 94, 362 (2016).
8. H. Mujahed, B. Nagy, Acta Crystallogr. Sect. A 72, 243 (2016).
9. M. Knor, R. Škrekovski, A. Tepeh, Discrete Appl. Math. 211, 121 (2016).
10. I. Gutman, N. Trinajstić, Chemical Physics Letters 17, 535 (1972).
11. B. Horoldagva, K. Das, T. Selenge, Discrete Appl. Math. 215, 146 (2016).
12. A. Ali, Z. Raza, A. Bhatti, Kuwait J. Sci. 43, 48 (2016).
13. S. Ediz, Mathematiche 71, 135 (2016).
14. M. Randić, Journal of the American Chemical Society 97, 6609 (1975).
15. R.K. Kincaid, S.J. Kunkler, M.D. Lamar, D.J. Phillips, Networks 67, 338 (2016).
16. A. Banerjee, R. Mehatari, Linear Algebra Appl. 505, 85 (2016).
17. R. Gu, F. Huang, X. Li, Trans. Comb. 5, 1 (2016).
18. B. Furtula, I. Gutman, J Math. Chem. 53, 1184 (2015).
19. W. Gao, M.K. Siddiqui, M. Imran, M.K. Jamil, M.R. Farahani, Saudi Pharmaceutical Journal 24, 258 (2016).
20. W. Gao, M.R. Farahani, L. Shi, Acta Medica Mediterranea, 32: 579 (2016).